\documentclass{IEEEtran4PSCC}
\usepackage{cite}
\ifCLASSINFOpdf
  \usepackage[pdftex]{graphicx}
\else
\fi
\usepackage[cmex10]{amsmath}
\usepackage{bm,mathtools}
\usepackage{url}
\usepackage[bookmarks=false]{hyperref}
\usepackage[draft]{minted}
\usepackage{float}
\floatstyle{ruled}
\newfloat{code}{thp}{lop}
\floatname{code}{Code Block}

\providecommand{\inlinecode}[1]{\texttt{#1}}

\usepackage{xcolor}

\begin{document}

%\title{PowerModels.jl: An Open-Source Toolkit \\ for Power Network Optimization}
\title{PowerModels.jl: An Open-Source Framework \\ for Exploring Power Flow Formulations}

%% To specify the authors when (number of affiliations <= 2)
\author{
\IEEEauthorblockN{Carleton Coffrin, Russell Bent, Kaarthik Sundar}
\IEEEauthorblockA{Los Alamos National Laboratory \\
Los Alamos, New Mexico, USA\\
\{cjc,rbent,kaarthik\}@lanl.gov}
\and
\IEEEauthorblockN{Yeesian Ng, Miles Lubin}
\IEEEauthorblockA{
%Operations Research Center\\ 
Massachusetts Institute of Technology\\
Cambridge, Massachusetts, USA\\
\{yeesian,mlubin\}@mit.edu}
}

%% To specify the authors when (number of affiliations > 2)
% \author{\IEEEauthorblockN{Author n.1\IEEEauthorrefmark{1},
% Author n.2\IEEEauthorrefmark{2},
% Author n.3\IEEEauthorrefmark{3}, 
% Author n.4\IEEEauthorrefmark{3} and
% Author n.5\IEEEauthorrefmark{4}}
% \IEEEauthorblockA{\IEEEauthorrefmark{1} Department Name of Organization A\\
% Name of the organization A,
% Address A\\ Emails if wanted}
% \IEEEauthorblockA{\IEEEauthorrefmark{2} Department Name of Organization B\\
% Name of the organization B,
% Address B\\ Emails if wanted}
% \IEEEauthorblockA{\IEEEauthorrefmark{3} Department Name of Organization C\\
% Name of the organization C,
% Address C\\ Emails if wanted}
% \IEEEauthorblockA{\IEEEauthorrefmark{4}Department Name of Organization D\\
% Name of the organization D,
% Address D\\ Emails if wanted}
% }

\maketitle

%\thanksto{bloop}

\begin{abstract}
In recent years, the power system research community has seen an explosion of novel methods for formulating and solving power network optimization problems.  These emerging methods range from new power flow approximations, which go beyond the traditional DC power flow by capturing reactive power, to convex relaxations, which provide solution quality and runtime performance guarantees.  Unfortunately, the sophistication of these emerging methods often presents a significant barrier to evaluating them on a wide variety of power system optimization applications.  To address this issue, this work proposes PowerModels, an open-source platform for comparing power flow formulations.  From its inception, PowerModels was designed to streamline the process of evaluating different power flow formulations on shared optimization problem specifications.  This work provides a brief introduction to the design of PowerModels, validates its implementation, and demonstrates its effectiveness with a proof-of-concept study analyzing five different formulations of the Optimal Power Flow problem.
\end{abstract}

%\vspace{-1cm}

% Use this to place sponsorships
%\thanksto{Applicable sponsors, if any, should be placed using the \emph{thanksto} command}

\section{Introduction \& Motivation}
In recent years, the power system research community has seen an explosion of novel methods for formulating and solving steady-state AC power network optimization problems.  
These emerging methods range from new power flow approximations, which go beyond the traditional DC power flow by capturing reactive power (e.g., LPAC \cite{LPAC_ijoc}, IV \cite{7299702}, LACR \cite{7447006}), to convex relaxations (e.g., Moment-Hierarchy \cite{moment_hierarchy}, SDP \cite{Bai2008383}, QC \cite{Hijazi2017}, SOC \cite{1664986}, CDF \cite{6102366}), which provide solution quality guarantees and leverage state-of-the-art convex optimization software.
Indeed, these emerging methods have demonstrated promising results on a wide range of problem domains, including Optimal Power Flow (OPF) \cite{7271127,5971792,moment_hierarchy}, Optimal Transmission Switching (OTS) \cite{Hijazi2017,pscc_ots}, Transmission Network Expansion Planning (TNEP) \cite{6407493,pscc_ep,7540988}, and Micro-Grid Design \cite{7947228,8023785}.
Unfortunately, a number of fundamental challenges have hindered ubiquitous access to these recent developments:
%\noindent
%\paragraph*{Challenges in Accessing Emerging Optimization Methods}
%
\begin{enumerate}
\item Many of the proposed methods are complex mathematical models including many auxiliary variables and constraints that are challenging to implement from scratch.  Furthermore, because of this complexity, the original authors' implementation is often required to precisely reproduce previous results.

\item In contrast to power system optimization algorithms built from scratch, these emerging methods often build on state-of-the-art optimization software, such as IPOPT \cite{Ipopt}, Mosek \cite{mosek}, or Gurobi \cite{gurobi}.  Typically these tools are tedious to use directly, so a modeling language, such as AMPL \cite{ampl}, GAMS \cite{gams}, OPL \cite{opl}, Pyomo \cite{hart2017pyomo}, YALMIP \cite{1393890}, or CVX \cite{cvx}, is often leveraged to make effective use of such software.

\item Matlab-based tools, such as {\sc Matpower} \cite{matpower}, have been widely successful as research and development baselines for power system analysis.  However, their reliance on Matlab can be a significant limitation.  Matlab licensing costs can be prohibitive in a nonacademic environment or when conducting large-scale experiments where hundreds of independent tasks need to be run in parallel.
\end{enumerate}

\paragraph*{The Advent of Julia \& JuMP}
Recently, Julia \cite{julia} has emerged as an open-source, high-level, high-performance programming language for numerical computing. Julia strives to have the ease of use of scripting languages like Matlab and Python while maintaining performance that is comparable to C.  In the context of this work, Julia is valuable because it presents a free and open-source alternative to Matlab.

JuMP \cite{DunningHuchetteLubin2017} is one of many open-source packages available in Julia and provides a modeling layer for optimization within Julia, similar to YALMIP and CVX in Matlab.  JuMP supports a wide variety of optimization problems, including linear programs (LP), mixed-integer programs (MIP), second-order cone programs (SOCP), semi-definite programs (SDP), nonlinear programs (NLP), and mixed-integer nonlinear programs (MINLP), which makes it an ideal modeling layer for the wide range of optimization problems that arise in power systems research.

\paragraph*{PowerModels}
Building on the success of Julia and JuMP, in this work we introduce PowerModels\footnote{\url{https://github.com/lanl-ansi/PowerModels.jl}}, a free and open-source toolkit for power network optimization with a focus on comparing power flow formulations.  At this time, PowerModels focuses on establishing a baseline implementation of steady-state power network optimization problems and includes implementations of Power Flow, OPF, OTS, and TNEP.  Each of these optimization problems can be considered with a variety of power flow formulations, including AC polar, AC rectangular, DC approximation, SOC relaxation, and the QC relaxation.  The correctness of the PowerModels implementation has been validated on hundreds of AC OPF benchmarks from both the PGLib-OPF \cite{pglib_opf} and NESTA \cite{nesta} benchmark archives.  The open-source nature of the Julia ecosystem makes replicating these results nearly effortless.  
%In time, we hope that PowerModels will be a valuable reference implementation for the wide variety of emerging optimization methods in the power system literature.  % covered in the conclusion
Furthermore, having all of these implementations in a common platform assists in fair and rigorous comparisons of different optimization methods from the power system literature.

The rest of the paper is organized as follows.  A brief review of mathematical programming is provided in Section \ref{sec:math_prog} to give context and motivation to this work.  Section \ref{sec:ps_math_prog} introduces the PowerModels framework, explains the design goals, and provides several examples of how it is used.  Section \ref{sec:study} validates the implementation of PowerModels to the established {\sc Matpower} software and demonstrates the effectiveness of the PowerModels framework with a proof-of-concept study comparing five different formulations of the OPF problem.  Section \ref{sec:conclusion} finishes with a few concluding remarks.

\begin{code}[t]
\caption{A JuMP v0.18 Model in Julia v0.6}
\label{code:jump_ex}
\begin{minted}{julia}
# Model
using JuMP
m = Model()

@variables m begin
  va[1:3]
  pg[1:3] >= 0
end

@objective(m, Min, 1*pg[1] + 10*pg[2] + 100*pg[3])

@constraints m begin
  va[1] == 0
  
  pg[1] - 1 == 10*(va[1]-va[2]) + 20*(va[1]-va[3])
  pg[2] - 2 == 10*(va[2]-va[1]) + 30*(va[2]-va[3])
  pg[3] - 4 == 20*(va[3]-va[1]) + 30*(va[3]-va[2])

  10*(va[1]-va[2]) <=  0.5
  10*(va[1]-va[2]) >= -0.5
  20*(va[1]-va[3]) <=  0.5
  20*(va[1]-va[3]) >= -0.5
end

# Solver
using Clp
setsolver(m, ClpSolver())
status = solve(m)
\end{minted}
\end{code}

\section{The Value of Mathematical Programming}
\label{sec:math_prog}

Since their inception in the 1970s, mathematical programming languages such as GAMS, AMPL, OPL, YALMIP, CVX, and JuMP have proven to be invaluable tools for specifying and solving a wide range of mathematical optimization problems.  One of the key features of these languages is the separation of the mathematical specification of the problem from the algorithmic task of solving it. It is important to note that, in such languages, it is easy to specify complex optimization problems that are impossible to solve in practice. However, a mathematically rigorous specification of such problems still has significant value in its own right \cite{go_comp,vn1553494}.

To make mathematical programming concrete, Code Block \ref{code:jump_ex} illustrates how a simple DC OPF problem can be modeled in the JuMP mathematical programming language. A key observation is that the complete mathematical specification is independent of the solution method, which is the CLP \cite{clp} linear programming solver in this case.  This abstraction layer is highly valuable for two reasons: (1) within the same problem class (e.g., LP), it is easy to explore a number of different solver technologies (e.g., CLP, GLPK, CPLEX, Gurobi);\footnote{See \url{http://www.juliaopt.org} for a list of solvers accessible via JuMP.} and (2) when modifications in the mathematical specification change the problem class (e.g., from an LP to an SDP), the only software change that is required is changing the solver technology (e.g., from the CLP linear solver to SCS \cite{ocpb:16} sdp solver).

\section{Power System Mathematical Programs}
\label{sec:ps_math_prog}

\begin{figure}[t]
    \begin{center}
    \includegraphics[width=8.8cm]{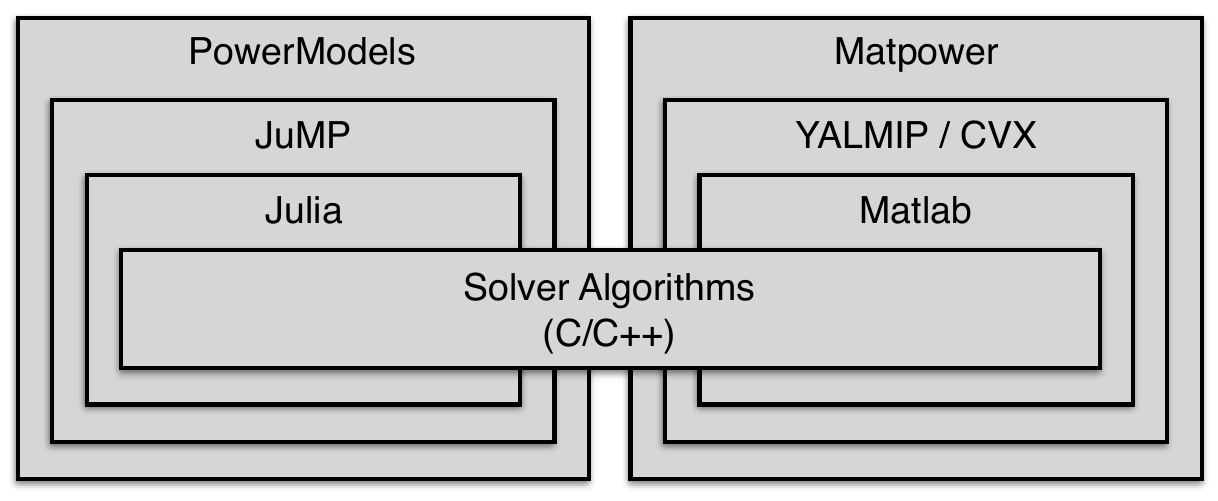}
    \end{center}
    \vspace{-0.4cm}
    \caption{Illustration of mathematical programming abstraction layers.}
    %\vspace{-0.2cm}
    \label{fig:layers}
\end{figure}

Inspired by the success of the abstractions developed by mathematical programming languages, a principal goal of PowerModels is to provide power system abstractions that will aid researchers in designing and comparing a wide range of power system optimization problems. As illustrated in Figure \ref{fig:layers}, PowerModels provides an abstraction layer on top of the JuMP mathematical programming language. This layer captures the structure of power systems mathematical programs, such as equations over complex numbers and network component objects such as buses, generators, and branches.

When developing any power system computational tool, two of the core decisions are (1) the types of system components to be supported, and (2) the mathematical model of those components. Based on the wide adoption of {\sc Matpower} in the research community, the first version of PowerModels elects the same component scope and mathematical model. Specifically, PowerModels supports (1) buses with one constant power load and one fixed shunt, (2) generators with polynomial cost functions, (3) $\pi$-equivalent branches and transformers, and (4) simple HVDC lines.  For a detailed description of each of these components, see the {\sc Matpower} documentation \cite{matpower}.  Selecting this formulation scope is also advantageous because of the significant amount of network data that has been curated in the {\sc Matpower} case format \cite{pglib_opf,nesta}.

In the pursuit of a power systems--specific mathematical programming framework, PowerModels makes two key observations that inform its design: (1) each power system optimization {\em problem} (e.g., OPF) has many {\em formulations} (e.g., AC in polar coordinates, DC approximation, or SOC relaxation), and (2) the AC power flow equations are most naturally defined in the space of complex numbers.  The following sections demonstrate how these observations manifest in a domain-specific mathematical programming framework for power system optimization.

\begin{code}[t]
\caption{PowerModels v0.5 Abstract OPF Model}
\label{code:opf}
\begin{minted}[numbersep=3pt,xleftmargin=8pt,linenos]{julia}
function post_opf(pm::GenericPowerModel)

  variable_voltage(pm)
  variable_generation(pm)
  variable_branch_flow(pm)
  variable_dcline_flow(pm)

  objective_min_fuel_cost(pm)

  constraint_voltage(pm)

  for i in ids(pm, :ref_buses)
    constraint_theta_ref(pm, i)
  end

  for i in ids(pm, :bus)
    constraint_kcl_shunt(pm, i)
  end

  for i in ids(pm, :branch)
    constraint_ohms_yt_from(pm, i)
    constraint_ohms_yt_to(pm, i)

    constraint_voltage_angle_difference(pm, i)

    constraint_thermal_limit_from(pm, i)
    constraint_thermal_limit_to(pm, i)
  end

  for i in ids(pm, :dcline)
    constraint_dcline(pm, i)
  end
end
\end{minted}
\end{code}

\subsection{Abstract Power System Problems}

Power system optimization problems in PowerModels are collections of functions with a syntax similar to that of JuMP mathematical models.  Code Block \ref{code:opf} provides the complete specification of the OPF problem, as prescribed by {\sc Matpower}.  This \inlinecode{post\_opf} function receives a \inlinecode{GenericPowerModel} object, which is analogous to a JuMP Model but is specialized to power systems.  The \inlinecode{GenericPowerModel} object includes various information about the PowerModel, most importantly the power network data.  This model object is then passed to a number of functions to build up the mathematical program.  Lines 3--6 define the decision variables.  Line 8 configures the objective function and Lines 10--32 add the problem constraints.  This model highlights two core features of power system mathematical models. First, all of the functions are defined over complex numbers.  For example, \inlinecode{variable\_voltage} initializes both the real and imaginary parts of the voltage variables and \inlinecode{constraint\_kcl\_shunt} captures Kirchoff's current law (KCL) on both active and reactive power. Second, the model constraints can be applied on a component-by-component basis; for example, KCL is applied for all buses whereas Ohm's law is applied for all branches.  This constraint organization highlights the network structure underlying mathematical programs in power systems.

To highlight the flexible nature of PowerModels as a modeling framework, Code Block \ref{code:ots} presents the OTS problem.  In this problem, formulation of an on/off branch indicator variable is added (i.e., \inlinecode{variable\_branch\_indicator}), and all of the branch flow constraints are extended to \inlinecode{\_on\_off} variants that incorporate this indicator variable.  Once those extensions are complete, the rest of the model leverages the same functions used in the OPF formulation (i.e., Code Block \ref{code:opf}).  The reuse of core abstractions in this example highlights the benefits of building a domain-specific mathematical programming framework for power systems.  

\begin{code}[t]
\caption{PowerModels v0.5 Abstract OTS Model}
\label{code:ots}
\begin{minted}{julia}
function post_ots(pm::GenericPowerModel)
  variable_branch_indicator(pm)
  variable_voltage_on_off(pm)
  variable_generation(pm)
  variable_branch_flow(pm)
  variable_dcline_flow(pm)

  objective_min_fuel_cost(pm)

  constraint_voltage_on_off(pm)

  for i in ids(pm, :ref_buses)
    constraint_theta_ref(pm, i)
  end

  for i in ids(pm, :bus)
    constraint_kcl_shunt(pm, i)
  end

  for i in ids(pm, :branch)
    constraint_ohms_yt_from_on_off(pm, i)
    constraint_ohms_yt_to_on_off(pm, i)

    constraint_voltage_angle_difference_on_off(pm, i)

    constraint_thermal_limit_from_on_off(pm, i)
    constraint_thermal_limit_to_on_off(pm, i)
  end

  for i in ids(pm, :dcline)
    constraint_dcline(pm, i)
  end
end
\end{minted}
\end{code}

% \begin{code}[t]
% \caption{PowerModels v0.5 Abstract TNEP Model}
% \label{code:tenp}
% \begin{minted}{julia}
% function post_tnep(pm::GenericPowerModel)
%   variable_branch_ne(pm)
%   variable_voltage(pm)
%   variable_voltage_ne(pm)
%   variable_generation(pm)
%   variable_branch_flow(pm)
%   variable_dcline_flow(pm)
%   variable_branch_flow_ne(pm)

%   objective_tnep_cost(pm)

%   constraint_voltage(pm)
%   constraint_voltage_ne(pm)
%   for i in ids(pm, :ref_buses)
%     constraint_theta_ref(pm, i)
%   end
%   for i in ids(pm, :bus)
%     constraint_kcl_shunt_ne(pm, i)
%   end
%   for i in ids(pm, :branch)
%     constraint_ohms_yt_from(pm, i)
%     constraint_ohms_yt_to(pm, i)
%     constraint_voltage_angle_difference(pm, i)
%     constraint_thermal_limit_from(pm, i)
%     constraint_thermal_limit_to(pm, i)
%   end
%   for i in ids(pm, :ne_branch)
%     constraint_ohms_yt_from_ne(pm, i)
%     constraint_ohms_yt_to_ne(pm, i)
%     constraint_voltage_angle_difference_ne(pm, i)
%     constraint_thermal_limit_from_ne(pm, i)
%     constraint_thermal_limit_to_ne(pm, i)
%   end
%   for i in ids(pm, :dcline)
%     constraint_dcline(pm, i)
%   end
% end
% \end{minted}
% \end{code}

\subsection{Power System Formulations}

The concept of power system {\em formulations} is introduced to transform generic problem specifications, such as the OPF and OTS problems from the previous section, into concrete mathematical programs, i.e., JuMP models. These formulations specify both the physics that will be used in the mathematical program (i.e., AC or DC) as well as the mathematical implementation of that physics (i.e., AC in polar form or AC in rectangular form).  PowerModels v0.5 includes the following formulations:
%\vspace{0.7cm}
\pagebreak

\begin{itemize}
\item {\em ACPPowerModel} - AC in polar coordinates \cite{ac_opf_origin}
\item {\em ACRPowerModel} - AC in rectangular coordinates \cite{6510541}
\item {\em ACTPowerModel} - AC in the w-theta space
\item {\em DCPPowerModel} - DC approximation \cite{matpower}
\item {\em DCPLLPowerModel} - DC approx. with line losses \cite{6345342}
\item {\em SOCWRPowerModel} - SOC relaxation \cite{1664986}
\item {\em QCWRPowerModel} - QC relaxation \cite{Hijazi2017}
\end{itemize}
All of these can be applied to any of the abstract problem formulations.

The key property that PowerModels assumes is that the combination of an abstract problem and a mathematical formulation results in a fully specified mathematical program that is then encoded as a JuMP model.  Code Block \ref{code:build_solve} illustrates how a problem and formulation are combined to build a JuMP model and then solve it.  The function \inlinecode{build\_generic\_model} combines the network data, a formulation, and a generic problem definition to build the JuMP model (i.e., \inlinecode{pm.model}).  The \inlinecode{solve\_generic\_model} function solves the JuMP model using the given solver ({\sc Ipopt} in this case) and then puts the raw numerical results in a data structure that is consistent with the PowerModels data format.  Because this build-and-solve process is fairly common, PowerModels provides the  \inlinecode{run\_opf} helper function to do one after the other.

\begin{code}[t]
\caption{Building and Running an AC-Polar OPF}
\label{code:build_solve}
\begin{minted}{julia}
using PowerModels
using Ipopt

# load the network data
case = PowerModels.parse_file("case24_ieee_rts.m")

# build the mathematical program
pm = build_generic_model(case, ACPPowerModel,
  PowerModels.post_opf)

# solve the mathematical program
result = solve_generic_model(pm, IpoptSolver())
\end{minted}
\end{code}

The intuition for how \inlinecode{build\_generic\_model} works is as follows.  A generic problem definition, such as OPF, is defined over a generic complex voltage product expression, i.e., $V_i V^*_j$.  Based on the given formulation, these generic expressions are replaced with a specific real number implementation.  For example, the {\em ACPPowerModel} formulation results in the following mapping:
\begin{align}
V_i V^*_j &\Rightarrow |V_i||V_j|\cos(\theta_i - \theta_j) + \bm i |V_i||V_j|\sin(\theta_i - \theta_j)
\end{align}
whereas the {\em ACRPowerModel} formulation has the mapping
\begin{align}
%V_i V^*_j &\Rightarrow (v^R_i + \bm i v^I_i)(v^R_j - \bm i v^R_i)
V_i V^*_j &\Rightarrow v^R_i v^R_j + v^I_i v^I_j + \bm i (v^I_i v^R_j - v^R_i v^I_j )
\end{align}
The {\em ACTPowerModel} formulation is an interesting case because it introduces constraints and new variables:
\begin{subequations}
\begin{align}
& V_i V^*_j \Rightarrow W^R_{ij} + \bm i W^I_{ij} \\
& (W^R_{ij})^2 + (W^I_{ij})^2 = W^R_{ii}W^R_{jj} \\
& W^I_{ij} = W^R_{ij} \tan(\theta_i - \theta_j)
\end{align}
\end{subequations}
and the traditional {\em DCPPowerModel} formulation is
\begin{align}
V_i V^*_j &\Rightarrow 0 + \bm i (\theta_i - \theta_j)
\end{align}
In general, this mapping may be more complex and add many auxiliary variables and constraints to the model, as is the case for the {\em QCWRPowerModel} formulation.  A detailed explanation of these notations can be found in \cite{7271127}.

\subsection{Comparing Formulations}

One of the interesting advantages of clearly separating the power system problem from the formulation is that it makes it easy to compare the effects of various formulations on a common problem of interest.  For example, consider Code Block \ref{code:run_forms}, in which the structure of each \inlinecode{run\_opf} call clearly indicates that the only thing changing is the formulation and not the underlying problem or solver.  This kind of structure mitigates implementation mistakes and assists in checking mathematical properties.  For example, the {\em SOCWRPowerModel} and {\em QCWRPowerModel} are designed to be relaxations of the nonconvex formulations {\em ACPPowerModel, ACRPowerModel}, and {\em ACTPowerModel}.  Hence, one expects that the objective value of both relaxations will be below the values of all the nonconvex formulations.

\begin{code}[t]
\caption{Running Various OPF Formulations}
\label{code:run_forms}
\begin{minted}{julia}
using PowerModels
using Ipopt

ipopt = IpoptSolver(tol=1e-6)

case = PowerModels.parse_file("case24_ieee_rts.m")

result_acp = run_opf(case, ACPPowerModel,   ipopt)
result_acr = run_opf(case, ACRPowerModel,   ipopt)
result_act = run_opf(case, ACTPowerModel,   ipopt)
result_soc = run_opf(case, SOCWRPowerModel, ipopt)
result_qc  = run_opf(case, QCWRPowerModel,  ipopt)
\end{minted}
\end{code}

\subsection{Proving Infeasibility}

Another interesting benefit of separating the power system optimization problem from the formulation is that it enables the use of convex relaxations as a diagnostic tool for proving that a problem's input data has no feasible solution, as noted in \cite{7271127}.  Such {\em proofs of infeasibility} are incredibly useful when debugging large network data sets.  The primary challenge is that large-scale power system optimization problems necessitate algorithms that provide only local optimality guarantees (e.g., Newton-Raphson and Interior Point Methods).  Therefore, if one solves an OPF problem with a nonconvex {\em ACRPowerModel} and the solver converges to an infeasible stationary point, it is not clear whether the problem is truly infeasible or whether the algorithm was unable to find a feasible solution.  In contrast, if one solves an OPF problem with a convex relaxation, such as {\em SOCWRPowerModel} or {\em QCWRPowerModel}, and the solver still converges to an infeasible stationary point, then the solver has proven that both the convex relaxation and the original nonconvex problem (e.g., {\em ACRPowerModel}) are infeasible.  The design of PowerModels makes it easy to perform this kind of problem data validation.

\subsection{User-Driven Extensions}

One of the advantages of building PowerModels in Julia is Julia's native package management system.  This allows users to develop extension packages on top of the PowerModels modeling layer and share those extensions with the community.  The PowerModelsAnnex\footnote{\url{https://github.com/lanl-ansi/PowerModelsAnnex.jl}} package provides an example of how this is accomplished and includes examples of how to extend PowerModels with new problem specifications and formulations.  Some preliminary user-driven activities include: PowerModelsReliability\footnote{\url{https://github.com/frederikgeth/PowerModelsReliability.jl}}, which includes problem extensions incorporating on-load tap changing transformers and load shedding; OPFRecourse\footnote{\url{https://github.com/lanl-ansi/OPFRecourse.jl}}, which includes problem extensions and algorithms for applying statistical learning methods to optimal power flow \cite{1801.07809,1802.09639}; And a worst-case $N$-$k$ contingency identification tool\footnote{\url{https://github.com/kaarthiksundar/n-k}} \cite{NET21806}.

\section{Validation and Proof-of-Concept Studies}
\label{sec:study}

The goals of this section are twofold. First, it describes a validation study to verify that PowerModels is comparable to the widely used {\sc Matpower} package. Second, it describes a proof-of-concept study to demonstrate the efficacy of using PowerModels for a comparison of different power flow formulations. Both studies consider the seminal OPF problem presented Code Block \ref{code:opf}.

\subsection{Test Cases and Computational Setting}

These studies consider the 108 power networks from the IEEE PES PGLib-OPF v17.08 benchmark library \cite{pglib_opf}.  The {\sc Matpower} v6.0 evaluation was conducted using Matlab R2017b and the default solvers provided for OPF problems.  The PowerModels v0.5 evaluation was conducted in Julia v0.6 using JuMP v0.18 \cite{DunningHuchetteLubin2017}.  In the interest of consistency, all of the PowerModels mathematical programs were solved with Ipopt \cite{Ipopt} using the HSL MA27 linear algebra solver \cite{hsl_lib} up to an optimality tolerance of $10^{-6}$.  All of the experiments were conducted on HPE ProLiant XL170r servers with two Intel CPUs @2.10 GHz and 128 GB of memory. The reported runtimes focus on the solver runtime and do not include the test case loading, model building, or Julia's JIT time. In total, these overheads add around 3--5 seconds. It was observed that on cases with $<$1000 buses, all OPF formulations solved in $<$1 second.  Hence, for brevity, these small cases are omitted.

\subsection{Validation Study}

Table \ref{tbl:mp_comp} presents a comparison of the AC-P and DC-P OPF solutions produced by {\sc Matpower} and PowerModels.  In the table,  {\em n.s.} indicates that no solution was found, {\em inf.} indicates that the solver proved that no solution exists, and -- indicates that the test case did not meet {\sc Matpower}'s data requirements.  Looking at the AC-P OPF results, in the cases where {\sc Matpower} finds a solution, its objective matches exactly with PowerModels, indicating that PowerModels correctly implements the AC power flow model.  Although {\sc Matpower}'s MIPS solver tends to be faster than Ipopt via PowerModels, it does appear to be less robust.  Looking at the DC-P OPF results, the objective values are similar but not exactly the same.  This is expected because the two tools implement slightly different variants of the DC power flow model.  Overall, these results suggest that PowerModels' OPF implementation is comparable to {\sc Matpower}'s.

\subsection{Formulation Study}

Table \ref{tbl:gaps_time} presents a comparison of the five formulations in Code Block \ref{code:run_forms}.  The AC-P model is used as the base case feasible solution to which all other formulations are compared.  For the nonconvex AC-R and AC-T formulations, the absolute difference in the objective value from AC-P is reported.  For the convex QC and SOC relaxations, the optimality gap is reported, that is,
%$100 * (\mbox{AC-P - Relaxation}) / \mbox{AC-P}$.
\begin{align}
100 * \frac{\mbox{AC-P - Relaxation}}{\mbox{AC-P}} \nonumber
\end{align}

The results are summarized as follows.
(1) Despite a lack of convergence guarantees, all three of the nonconvex models AC-P, AC-R, and AC-T converge to very similar solutions.  Of these three, AC-R tends to be the fastest; this is consistent with the results of \cite{acopf_comp}.
(2) In terms of the relaxations, the results confirm that the QC is a stronger relaxation than the SOC, as shown in \cite{7271127}, but there is a clear runtime benefit for the simpler SOC relaxation.
(3) The consistency of both the nonconvex models and convex relaxations suggests that the PowerModels implementation of these models is correct.

% \subsection{Replication of the Results}
% [TODO]
% Bloop

%%%%%%%
%%%%%%%
% Generated by opf-mp-tbl.py in https://github.lanlytics.com/cjc/pms-exper-pub
% pglib-mp-comp.tex
%%%%%%%
%%%%%%%
\begin{table*}[t]
\small
\caption{Quality and Runtime Comparison of Matpower and PowerModels on AC and DC Optimal Power Flow}
\begin{tabular}{|r|r|r||r|r||r|r||r|r|r|r|r|r|r|r|r|r|r|r|r|r|}
\hline
 & & & \multicolumn{4}{c||}{\$/h} & \multicolumn{4}{c|}{Runtime (seconds)} \\
 & & & MP & PM & MP & PM & MP & PM & MP & PM \\
Test Case & $|N|$ & $|E|$ & AC-P & AC-P & DC & DC & AC-P & AC-P & DC & DC \\
\hline
\hline
\multicolumn{11}{|c|}{Typical Operating Conditions (TYP)} \\
\hline
case1354\_pegase & 1354 & 1991 & 1.3640e+06 & 1.3640e+06 & 1.3141e+06 & 1.3140e+06 & 4 & 6 & 2 & $<$1 \\
\hline
case1888\_rte & 1888 & 2531 & n.s. & 1.5654e+06 & 1.5111e+06 & 1.5111e+06 & 6 & 14 & 2 & $<$1 \\
\hline
case1951\_rte & 1951 & 2596 & n.s. & 2.3753e+06 & 2.3128e+06 & 2.3128e+06 & 3 & 18 & 2 & $<$1 \\
\hline
case2383wp\_k & 2383 & 2896 & 1.8685e+06 & 1.8685e+06 & 1.7968e+06 & 1.8041e+06 & 5 & 9 & 2 & $<$1 \\
\hline
case2736sp\_k & 2736 & 3504 & 1.3079e+06 & 1.3079e+06 & 1.2760e+06 & 1.2760e+06 & 4 & 8 & 2 & $<$1 \\
\hline
case2737sop\_k & 2737 & 3506 & 7.7763e+05 & 7.7763e+05 & 7.6401e+05 & 7.6401e+05 & 5 & 6 & 2 & $<$1 \\
\hline
case2746wop\_k & 2746 & 3514 & 1.2083e+06 & 1.2083e+06 & 1.1782e+06 & 1.1782e+06 & 5 & 7 & 2 & $<$1 \\
\hline
case2746wp\_k & 2746 & 3514 & 1.6318e+06 & 1.6318e+06 & 1.5814e+06 & 1.5814e+06 & 5 & 7 & 3 & $<$1 \\
\hline
case2848\_rte & 2848 & 3776 & n.s. & 1.3847e+06 & 1.3636e+06 & 1.3636e+06 & 8 & 19 & 2 & $<$1 \\
\hline
case2868\_rte & 2868 & 3808 & n.s. & 2.2599e+06 & 2.2053e+06 & 2.2053e+06 & 12 & 21 & 3 & $<$1 \\
\hline
case2869\_pegase & 2869 & 4582 & 2.6050e+06 & 2.6050e+06 & 2.5167e+06 & 2.5166e+06 & 8 & 14 & 3 & $<$1 \\
\hline
case3012wp\_k & 3012 & 3572 & 2.6008e+06 & 2.6008e+06 & 2.5143e+06 & 2.5090e+06 & 7 & 11 & 3 & $<$1 \\
\hline
case3120sp\_k & 3120 & 3693 & 2.1457e+06 & 2.1457e+06 & 2.0891e+06 & 2.0880e+06 & 7 & 11 & 3 & $<$1 \\
\hline
case3375wp\_k & 3375 & 4161 & n.s. & 7.4357e+06 & 7.3166e+06 & 7.3170e+06 & 2 & 14 & 3 & $<$1 \\
\hline
case6468\_rte & 6468 & 9000 & n.s. & 2.2623e+06 & 2.1797e+06 & 2.1619e+06 & 7 & 80 & 4 & $<$1 \\
\hline
case6470\_rte & 6470 & 9005 & n.s. & 2.5558e+06 & 2.4520e+06 & 2.4454e+06 & 32 & 47 & 4 & 2 \\
\hline
case6495\_rte & 6495 & 9019 & n.s. & 3.4777e+06 & 3.0085e+06 & 2.8481e+06 & 16 & 89 & 4 & 2 \\
\hline
case6515\_rte & 6515 & 9037 & n.s. & 3.1971e+06 & 2.9493e+06 & 2.8484e+06 & 30 & 73 & 4 & 2 \\
\hline
case9241\_pegase & 9241 & 16049 & n.s. & 6.7747e+06 & 6.5411e+06 & 6.5179e+06 & 68 & 62 & 7 & 2 \\
\hline
case13659\_pegase & 13659 & 20467 & 1.0781e+07 & 1.0781e+07 & 1.0587e+07 & 1.0565e+07 & 51 & 96 & 8 & 3 \\
\hline
\hline
\multicolumn{11}{|c|}{Congested Operating Conditions (API)} \\
\hline
case1354\_pegase\_\_api & 1354 & 1991 & -- & 1.8041e+06 & 1.7479e+06 & 1.7485e+06 & -- & 6 & 2 & $<$1 \\
\hline
case1888\_rte\_\_api & 1888 & 2531 & n.s. & 2.2566e+06 & 2.1930e+06 & 2.1930e+06 & 9 & 10 & 2 & $<$1 \\
\hline
case1951\_rte\_\_api & 1951 & 2596 & -- & 2.8005e+06 & n.s. & 2.7027e+06 & -- & 124 & 2 & $<$1 \\
\hline
case2383wp\_k\_\_api & 2383 & 2896 & 2.7913e+05 & 2.7913e+05 & 2.7913e+05 & 2.7913e+05 & 3 & 31 & 2 & $<$1 \\
\hline
case2736sp\_k\_\_api & 2736 & 3504 & 6.3847e+05 & 6.3847e+05 & 6.1252e+05 & 5.9774e+05 & 5 & 9 & 3 & $<$1 \\
\hline
case2737sop\_k\_\_api & 2737 & 3506 & 4.0282e+05 & 4.0282e+05 & 3.7767e+05 & 3.7570e+05 & 5 & 8 & 3 & $<$1 \\
\hline
case2746wop\_k\_\_api & 2746 & 3514 & 5.1166e+05 & 5.1166e+05 & 5.1166e+05 & 5.1166e+05 & 3 & 3 & 3 & $<$1 \\
\hline
case2746wp\_k\_\_api & 2746 & 3514 & 5.8183e+05 & 5.8183e+05 & 5.8183e+05 & 5.8183e+05 & 4 & 5 & 3 & $<$1 \\
\hline
case2848\_rte\_\_api & 2848 & 3776 & -- & 1.7169e+06 & 1.6822e+06 & 1.6822e+06 & -- & 35 & 3 & $<$1 \\
\hline
case2868\_rte\_\_api & 2868 & 3808 & -- & 2.7159e+06 & n.s. & 2.6357e+06 & -- & 29 & 3 & $<$1 \\
\hline
case2869\_pegase\_\_api & 2869 & 4582 & -- & 3.3185e+06 & 3.2137e+06 & 3.2154e+06 & -- & 16 & 3 & $<$1 \\
\hline
case3012wp\_k\_\_api & 3012 & 3572 & 7.2887e+05 & 7.2887e+05 & 7.2887e+05 & 7.2887e+05 & 7 & 5 & 3 & $<$1 \\
\hline
case3120sp\_k\_\_api & 3120 & 3693 & 9.2026e+05 & 9.2026e+05 & 8.9103e+05 & 8.5997e+05 & 7 & 15 & 3 & $<$1 \\
\hline
case3375wp\_k\_\_api & 3375 & 4161 & n.s. & 5.8861e+06 & 5.8117e+06 & 5.7641e+06 & 2 & 14 & 3 & $<$1 \\
\hline
case6468\_rte\_\_api & 6468 & 9000 & n.s. & 2.7102e+06 & 2.6078e+06 & 2.6081e+06 & 45 & 84 & 4 & $<$1 \\
\hline
case6470\_rte\_\_api & 6470 & 9005 & n.s. & 3.1603e+06 & 3.0266e+06 & 3.0333e+06 & 18 & 58 & 4 & $<$1 \\
\hline
case6495\_rte\_\_api & 6495 & 9019 & n.s. & 3.6263e+06 & 3.4265e+06 & 3.4186e+06 & 7 & 77 & 4 & 2 \\
\hline
case6515\_rte\_\_api & 6515 & 9037 & n.s. & 3.5904e+06 & 3.3662e+06 & 3.3777e+06 & 42 & 72 & 4 & 2 \\
\hline
case9241\_pegase\_\_api & 9241 & 16049 & -- & 8.2656e+06 & 7.9822e+06 & 7.9822e+06 & -- & 73 & 13 & 2 \\
\hline
case13659\_pegase\_\_api & 13659 & 20467 & -- & 1.1209e+07 & 1.0903e+07 & 1.0895e+07 & -- & 79 & 11 & 3 \\
\hline
\hline
\multicolumn{11}{|c|}{Small Angle Difference Conditions (SAD)} \\
\hline
case1354\_pegase\_\_sad & 1354 & 1991 & 1.3646e+06 & 1.3646e+06 & n.s. & inf. & 4 & 6 & 2 & $<$1 \\
\hline
case1888\_rte\_\_sad & 1888 & 2531 & n.s. & 1.5806e+06 & 1.5122e+06 & 1.5123e+06 & 7 & 16 & 2 & $<$1 \\
\hline
case1951\_rte\_\_sad & 1951 & 2596 & n.s. & 2.3820e+06 & n.s. & inf. & 3 & 25 & 2 & $<$1 \\
\hline
case2383wp\_k\_\_sad & 2383 & 2896 & 1.9165e+06 & 1.9165e+06 & n.s. & inf. & 5 & 10 & 2 & 2 \\
\hline
case2736sp\_k\_\_sad & 2736 & 3504 & 1.3294e+06 & 1.3294e+06 & n.s. & inf. & 5 & 10 & 2 & 2 \\
\hline
case2737sop\_k\_\_sad & 2737 & 3506 & 7.9267e+05 & 7.9267e+05 & n.s. & inf. & 5 & 9 & 2 & 2 \\
\hline
case2746wop\_k\_\_sad & 2746 & 3514 & 1.2344e+06 & 1.2344e+06 & n.s. & inf. & 5 & 8 & 2 & 2 \\
\hline
case2746wp\_k\_\_sad & 2746 & 3514 & 1.6674e+06 & 1.6674e+06 & n.s. & inf. & 5 & 9 & 3 & 19 \\
\hline
case2848\_rte\_\_sad & 2848 & 3776 & n.s. & 1.3879e+06 & n.s. & inf. & 6 & 21 & 3 & 2 \\
\hline
case2868\_rte\_\_sad & 2868 & 3808 & n.s. & 2.2707e+06 & n.s. & inf. & 7 & 20 & 3 & 2 \\
\hline
case2869\_pegase\_\_sad & 2869 & 4582 & 2.6198e+06 & 2.6198e+06 & n.s. & inf. & 9 & 14 & 3 & 4 \\
\hline
case3012wp\_k\_\_sad & 3012 & 3572 & 2.6213e+06 & 2.6213e+06 & n.s. & inf. & 7 & 12 & 3 & 3 \\
\hline
case3120sp\_k\_\_sad & 3120 & 3693 & 2.1755e+06 & 2.1755e+06 & n.s. & inf. & 7 & 14 & 3 & 3 \\
\hline
case3375wp\_k\_\_sad & 3375 & 4161 & n.s. & 7.4357e+06 & 7.3181e+06 & 7.3197e+06 & 2 & 14 & 3 & $<$1 \\
\hline
case6468\_rte\_\_sad & 6468 & 9000 & n.s. & 2.2623e+06 & 2.1797e+06 & 2.1619e+06 & 7 & 131 & 4 & $<$1 \\
\hline
case6470\_rte\_\_sad & 6470 & 9005 & n.s. & 2.5597e+06 & 2.4595e+06 & 2.4483e+06 & 22 & 48 & 4 & 2 \\
\hline
case6495\_rte\_\_sad & 6495 & 9019 & n.s. & 3.4777e+06 & 3.0995e+06 & 2.8482e+06 & 12 & 88 & 4 & 2 \\
\hline
case6515\_rte\_\_sad & 6515 & 9037 & n.s. & 3.2679e+06 & 3.1394e+06 & 2.8486e+06 & 11 & 78 & 4 & 2 \\
\hline
case9241\_pegase\_\_sad & 9241 & 16049 & n.s. & 6.9170e+06 & n.s. & inf. & 38 & 70 & 5 & 10 \\
\hline
case13659\_pegase\_\_sad & 13659 & 20467 & 1.0901e+07 & 1.0901e+07 & n.s. & inf. & 51 & 79 & 6 & 12 \\
\hline
\end{tabular}\\
\label{tbl:mp_comp}
\end{table*}

%%%%%%%
%%%%%%%
% Generated by opf-delta-tbl.py in https://github.lanlytics.com/cjc/pms-exper-pub
% pglib-form-comp.tex
%%%%%%%
%%%%%%%
\begin{table*}[t]
\small
\caption{Quality and Runtime Results for AC Optimal Power Flow Formulations and Convex Relaxations}
\begin{tabular}{|r|r|r||r||r|r||r|r||r|r|r|r|r|r|r|r|r|r|r|r|r|}
\hline
& & & \$/h & \multicolumn{2}{c||}{$\Delta$\$/h} & \multicolumn{2}{c||}{Gap (\%)} & \multicolumn{5}{c|}{Runtime (seconds)} \\
Test Case & $|N|$ & $|E|$ & AC-P & AC-R & AC-T & QC & SOC & AC-P & AC-R & AC-T & QC & SOC \\
\hline
\hline
\multicolumn{13}{|c|}{Typical Operating Conditions (TYP)} \\
\hline
case1354\_pegase & 1354 & 1991 & 1.3640e+06 & -0.00 & 0.00 & 2.40 & 2.41 & 6 & 2 & 3 & 6 & 3 \\
\hline
case1888\_rte & 1888 & 2531 & 1.5654e+06 & 0.07 & 0.07 & 1.82 & 1.82 & 14 & 31 & 27 & 9 & 5 \\
\hline
case1951\_rte & 1951 & 2596 & 2.3753e+06 & 0.00 & 0.00 & 0.12 & 0.13 & 18 & 7 & 21 & 10 & 6 \\
\hline
case2383wp\_k & 2383 & 2896 & 1.8685e+06 & -0.00 & 0.00 & 0.99 & 1.05 & 9 & 4 & 127 & 10 & 6 \\
\hline
case2736sp\_k & 2736 & 3504 & 1.3079e+06 & 0.00 & 0.00 & 0.29 & 0.30 & 8 & 3 & 19 & 10 & 4 \\
\hline
case2737sop\_k & 2737 & 3506 & 7.7763e+05 & 0.00 & 0.00 & 0.25 & 0.26 & 6 & 3 & 27 & 9 & 4 \\
\hline
case2746wop\_k & 2746 & 3514 & 1.2083e+06 & 0.00 & 0.00 & 0.36 & 0.37 & 7 & 3 & 22 & 10 & 43 \\
\hline
case2746wp\_k & 2746 & 3514 & 1.6318e+06 & 0.00 & 0.00 & 0.32 & 0.33 & 7 & 3 & 24 & 10 & 5 \\
\hline
case2848\_rte & 2848 & 3776 & 1.3847e+06 & -0.00 & 0.00 & 0.12 & 0.12 & 19 & 6 & 9 & 14 & 7 \\
\hline
case2868\_rte & 2868 & 3808 & 2.2599e+06 & -0.00 & 0.00 & 0.11 & 0.11 & 21 & 7 & 11 & 14 & 8 \\
\hline
case2869\_pegase & 2869 & 4582 & 2.6050e+06 & -0.00 & 0.00 & 1.07 & 1.08 & 14 & 6 & 8 & 20 & 7 \\
\hline
case3012wp\_k & 3012 & 3572 & 2.6008e+06 & 0.00 & 0.00 & 0.98 & 1.03 & 11 & 5 & 16 & 14 & 8 \\
\hline
case3120sp\_k & 3120 & 3693 & 2.1457e+06 & 0.00 & 0.00 & 0.54 & 0.55 & 11 & 5 & 20 & 15 & 6 \\
\hline
case3375wp\_k & 3375 & 4161 & 7.4357e+06 & 0.00 & err. & 0.50 & 0.52 & 14 & 16 & 41 & 44 & 26 \\
\hline
case6468\_rte & 6468 & 9000 & 2.2623e+06 & 0.00 & 0.00 & 1.07 & 1.07 & 80 & 49 & 75 & 69 & 27 \\
\hline
case6470\_rte & 6470 & 9005 & 2.5558e+06 & 0.00 & 0.00 & 1.95 & 1.96 & 47 & 19 & 124 & 48 & 23 \\
\hline
case6495\_rte & 6495 & 9019 & 3.4777e+06 & 0.03 & 0.09 & 16.73 & 16.75 & 89 & 35 & 181 & 52 & 24 \\
\hline
case6515\_rte & 6515 & 9037 & 3.1971e+06 & 0.02 & 0.04 & 7.86 & 7.87 & 73 & 28 & 249 & 56 & 22 \\
\hline
case9241\_pegase & 9241 & 16049 & 6.7747e+06 & 0.00 & 0.00 & 1.99 & 2.84 & 62 & 24 & 525 & 121 & 37 \\
\hline
case13659\_pegase & 13659 & 20467 & 1.0781e+07 & 0.00 & 0.00 & 0.95 & 1.35 & 96 & 67 & 499 & 131 & 45 \\
\hline
\hline
\multicolumn{13}{|c|}{Congested Operating Conditions (API)} \\
\hline
case1354\_pegase\_\_api & 1354 & 1991 & 1.8041e+06 & -0.00 & 0.00 & 0.70 & 0.71 & 6 & 2 & 3 & 6 & 3 \\
\hline
case1888\_rte\_\_api & 1888 & 2531 & 2.2566e+06 & 0.02 & 0.03 & 0.47 & 0.47 & 10 & 4 & 11 & 14 & 5 \\
\hline
case1951\_rte\_\_api & 1951 & 2596 & 2.8005e+06 & -0.03 & 0.01 & 0.60 & 0.62 & 124 & 4 & 10 & 10 & 5 \\
\hline
case2383wp\_k\_\_api & 2383 & 2896 & 2.7913e+05 & 0.00 & -0.00 & 0.01 & 0.01 & 31 & 2 & 2 & 3 & $<$1 \\
\hline
case2736sp\_k\_\_api & 2736 & 3504 & 6.3847e+05 & 0.00 & 0.00 & 12.82 & 12.83 & 9 & 4 & 5 & 10 & 3 \\
\hline
case2737sop\_k\_\_api & 2737 & 3506 & 4.0282e+05 & 0.00 & 0.00 & 10.00 & 10.01 & 8 & 3 & 51 & 9 & 3 \\
\hline
case2746wop\_k\_\_api & 2746 & 3514 & 5.1166e+05 & -0.00 & -0.00 & 0.01 & 0.01 & 3 & 2 & 2 & 4 & 2 \\
\hline
case2746wp\_k\_\_api & 2746 & 3514 & 5.8183e+05 & -0.00 & -0.00 & 0.01 & 0.00 & 5 & 2 & 3 & 5 & 2 \\
\hline
case2848\_rte\_\_api & 2848 & 3776 & 1.7169e+06 & -0.00 & 0.00 & 0.18 & 0.18 & 35 & 11 & 40 & 15 & 7 \\
\hline
case2868\_rte\_\_api & 2868 & 3808 & 2.7159e+06 & -0.00 & 0.00 & 0.19 & 0.20 & 29 & 10 & 14 & 16 & 6 \\
\hline
case2869\_pegase\_\_api & 2869 & 4582 & 3.3185e+06 & -0.00 & 0.00 & 0.82 & 0.84 & 16 & 7 & 11 & 20 & 8 \\
\hline
case3012wp\_k\_\_api & 3012 & 3572 & 7.2887e+05 & 0.00 & -0.00 & 0.00 & 0.00 & 5 & 2 & 3 & 6 & 2 \\
\hline
case3120sp\_k\_\_api & 3120 & 3693 & 9.2026e+05 & 0.00 & 0.00 & 24.92 & 24.95 & 15 & 6 & 9 & 14 & 4 \\
\hline
case3375wp\_k\_\_api & 3375 & 4161 & 5.8861e+06 & 0.05 & 0.00 & 9.46 & 9.55 & 14 & 38 & 259 & 21 & 14 \\
\hline
case6468\_rte\_\_api & 6468 & 9000 & 2.7102e+06 & -0.00 & 0.00 & 0.41 & 0.42 & 84 & 40 & 84 & 54 & 29 \\
\hline
case6470\_rte\_\_api & 6470 & 9005 & 3.1603e+06 & -0.00 & 0.00 & 0.82 & 0.84 & 58 & 24 & 92 & 40 & 23 \\
\hline
case6495\_rte\_\_api & 6495 & 9019 & 3.6263e+06 & 0.31 & 0.31 & 3.24 & 3.28 & 77 & 27 & 134 & 49 & 23 \\
\hline
case6515\_rte\_\_api & 6515 & 9037 & 3.5904e+06 & 0.22 & 0.23 & 2.53 & 2.56 & 72 & 30 & 174 & 53 & 22 \\
\hline
case9241\_pegase\_\_api & 9241 & 16049 & 8.2656e+06 & 0.99 & 0.00 & 1.70 & 2.59 & 73 & 32 & 1536 & 141 & 39 \\
\hline
case13659\_pegase\_\_api & 13659 & 20467 & 1.1209e+07 & 0.00 & 0.00 & 1.21 & 1.91 & 79 & 37 & 768 & 122 & 56 \\
\hline
\hline
\multicolumn{13}{|c|}{Small Angle Difference Conditions (SAD)} \\
\hline
case1354\_pegase\_\_sad & 1354 & 1991 & 1.3646e+06 & -0.00 & 0.00 & 2.37 & 2.45 & 6 & 2 & 3 & 6 & 3 \\
\hline
case1888\_rte\_\_sad & 1888 & 2531 & 1.5806e+06 & 0.09 & 0.09 & 2.73 & 2.74 & 16 & 8 & 10 & 8 & 5 \\
\hline
case1951\_rte\_\_sad & 1951 & 2596 & 2.3820e+06 & 0.01 & 0.00 & 0.37 & 0.40 & 25 & 6 & 14 & 9 & 5 \\
\hline
case2383wp\_k\_\_sad & 2383 & 2896 & 1.9165e+06 & 0.01 & 0.00 & 2.16 & 3.13 & 10 & 4 & 60 & 10 & 6 \\
\hline
case2736sp\_k\_\_sad & 2736 & 3504 & 1.3294e+06 & 0.00 & 0.00 & 1.53 & 1.80 & 10 & 4 & 23 & 11 & 5 \\
\hline
case2737sop\_k\_\_sad & 2737 & 3506 & 7.9267e+05 & 0.00 & 0.00 & 1.92 & 2.10 & 9 & 4 & 14 & 9 & 4 \\
\hline
case2746wop\_k\_\_sad & 2746 & 3514 & 1.2344e+06 & 0.00 & 0.00 & 2.00 & 2.37 & 8 & 4 & 20 & 8 & 4 \\
\hline
case2746wp\_k\_\_sad & 2746 & 3514 & 1.6674e+06 & 0.00 & 0.00 & 1.68 & 2.21 & 9 & 4 & 18 & 10 & 5 \\
\hline
case2848\_rte\_\_sad & 2848 & 3776 & 1.3879e+06 & 0.00 & 0.00 & 0.27 & 0.29 & 21 & 8 & 13 & 14 & 5 \\
\hline
case2868\_rte\_\_sad & 2868 & 3808 & 2.2707e+06 & 0.00 & 0.00 & 0.50 & 0.53 & 20 & 8 & 10 & 15 & 6 \\
\hline
case2869\_pegase\_\_sad & 2869 & 4582 & 2.6198e+06 & 0.00 & 0.00 & 1.39 & 1.49 & 14 & 6 & 8 & 31 & 7 \\
\hline
case3012wp\_k\_\_sad & 3012 & 3572 & 2.6213e+06 & 0.01 & 0.00 & 1.41 & 1.62 & 12 & 6 & 19 & 15 & 6 \\
\hline
case3120sp\_k\_\_sad & 3120 & 3693 & 2.1755e+06 & 0.01 & 0.00 & 1.42 & 1.61 & 14 & 6 & 27 & 15 & 6 \\
\hline
case3375wp\_k\_\_sad & 3375 & 4161 & 7.4357e+06 & 0.00 & 0.00 & 0.47 & 0.52 & 14 & 18 & 23 & 27 & 19 \\
\hline
case6468\_rte\_\_sad & 6468 & 9000 & 2.2623e+06 & 0.00 & 0.00 & 1.05 & 1.06 & 131 & 41 & 138 & 67 & 28 \\
\hline
case6470\_rte\_\_sad & 6470 & 9005 & 2.5597e+06 & 0.00 & 0.00 & 2.03 & 2.08 & 48 & 19 & 65 & 46 & 22 \\
\hline
case6495\_rte\_\_sad & 6495 & 9019 & 3.4777e+06 & 0.03 & 0.09 & 16.63 & 16.75 & 88 & 34 & 160 & 54 & 23 \\
\hline
case6515\_rte\_\_sad & 6515 & 9037 & 3.2679e+06 & 0.02 & 0.03 & 9.82 & 9.87 & 78 & 30 & 144 & 53 & 22 \\
\hline
case9241\_pegase\_\_sad & 9241 & 16049 & 6.9170e+06 & 27.87 & 27.87 & 3.49 & 3.56 & 70 & 29 & 393 & 143 & 37 \\
\hline
case13659\_pegase\_\_sad & 13659 & 20467 & 1.0901e+07 & 0.01 & 0.01 & 1.70 & 1.74 & 79 & 126 & 635 & 135 & 910 \\
\hline
\end{tabular}\\
\label{tbl:gaps_time}
\end{table*}

\section{Conclusion}
\label{sec:conclusion}

In this work, PowerModels is proposed as a domain-specific mathematical programming framework for power system optimization.  A core insight of this framework is that it is possible to factor power system mathematical programs into two independent components, problem specifications and mathematical formulations.  The success of this approach was demonstrated by a proof-of-concept comparison of five established formulations of the OPF problem.  In the future we hope that PowerModels will mature well beyond the current version and will include a wide selection of power flow formulations from the literature, especially the well-known SDP and Moment-Hierarchy relaxations.  We encourage the power  systems research community to share their experiences using PowerModels and to contribute novel problem specifications and mathematical formulations.

\clearpage





% \bibitem{Shell}
% M.~Shell, \emph{How to Use the IEEEtran Latex Class}, Latex Archive Contents, \verb+http://www.ieee.org/conferences_events/+ \verb+conferences/publishing/templates.htm+

% \end{thebibliography}

\noindent
LA-UR-17-29326

\end{document}